\documentclass[a4paper,12pt]{article}
\usepackage{amssymb}
 \usepackage{amsthm}
\usepackage{latexsym}
\usepackage{amsmath}
\usepackage{amssymb}
\newcommand{\udots}{\mathinner{\mskip1mu\raise1pt\vbox{\kern7pt\hbox{.}}
\mskip2mu\raise4pt\hbox{.}\mskip2mu\raise7pt\hbox{.}\mskip1mu}}
\usepackage{amsmath}

\newtheorem{theorem}{Theorem}[section]
\newtheorem{corollary}{Corollary}[section]
\newtheorem{lemma}{Lemma}[section]
\newtheorem{remark}{Remark}[section]

\begin{document}

\baselineskip=21pt

\title{\bf \Large {The mixability of elliptical distributions and   log-elliptical distributions}}
{{\author{\normalsize{Xiaoqian Zhang\;\;\; Xiang Li\;\;\; Chuancun Yin}\\
{\normalsize\it  School of Statistics,  Qufu Normal University}\\
\noindent{\normalsize\it Shandong 273165, China}\\
e-mail:  ccyin@qfnu.edu.cn}}}
\maketitle
\vskip0.01cm
%\newpage
\noindent{{\bf ABSTRACT}  The concept of $\phi$-complete mixability and $\phi$-joint mixability was first introduced  in Bignozzi and Puccetti (2015), which is an extension of complete and joint  mixability. Following  Bignozzi and Puccetti (2015), we consider two     more general  cases of $\phi$  and investigate the $\phi$-joint mixability  for   elliptical distributions and logarithmic elliptical distributions. Some sufficient conditions for the  $\phi$-joint mixability of some  distributions are investigated.
 In addition, a conjecture on the  uniqueness of the center  of  $\phi$-joint  mixability  and the forms of the densities for some  elliptical distributions are given.}



%\newpage

%\noindent{\bf 1.~~Introduction}
\numberwithin{equation}{section}
\section{Introduction}\label{intro}

%The concept of  complete mixability and  joint mixability has been studied by many researchers in recent years because of its important application in many relevant fields.
 The  concept of complete mixability for a univariate distribution was first introduced by Wang and Wang (2011) and then extended to the notion of joint mixability of an
arbitrary set of distributions by  Wang,  Peng and  Yang (2013).
Suppose $n$ is a positive integer, distribution functions $F_1,\cdots, F_n$ on $\mathbb{R}$ are  said to be jointly mixable with index $n$ if there exist $n$ random variables $X_1,\cdots, X_n$  such that $X_i\sim F_i$, $1\le i\le n$, and
$P(X_1 + \cdots + X_n =C) = 1$
for some $C\in \Bbb{R}$. If $F_i=F$, $1\le i\le n$, then $F$ is said to be $n$-completely mixable and
 $(X_1,\cdots,X_n)$ a complete mix.  Any such $C$ is called a joint center of   $(F_1,\cdots, F_n)$.
 The concept of complete mixability is related to
some optimization problems in the theory of optimal couplings. For more details on the problems and a brief history of the concept of the
 mixability, we refer to the  papers of  Puccetti,  Wang and Wang (2012), Wang (2015) and Wang and Wang (2016). Bignozzi and Puccetti (2015) extended the concept of joint mixability  and introduce the concept of $\phi$-joint mixability as follows.
  Let $\phi:\mathbb{R}^n{\rightarrow}\mathbb{R}$ be a measurable function. An $n$-tuple of univariate distribution functions ($F_1, \cdots, F_n$) is said to be $\phi$-jointly mixable with index $n$ if there exist $n$ random variables $X_1,\cdots,X_n$ such that $X_i\sim\,F_i,$ $1\,{\leq}\,i\,{\leq}\,n$, and
\begin{eqnarray}
P(\phi(X_1,\cdots,X_n)=C)=1
\end{eqnarray}
for some $C\in\mathbb{R}.$ Any such $C$ is called a center of the $\phi$-jointly mixable distributions  and   any random vector $(X_1, \cdots, X_n)$
satisfying (1.1) is called a joint $\phi$-mix. If $F_i=F$, $1\le i\le n$, then $F$ is said to be $\phi$-completely mixable with index $n$  and  the random vector $(X_1, \cdots, X_n)$ is called a complete $\phi$-mix.

 In this paper, we consider the following two kind functions:
\begin{eqnarray}
\phi_1(x_1,\cdots,x_n)=h({\alpha_1}x_1+\cdots+{\alpha_n}x_n),\ \alpha_i>0,\ x_i\,\in\,\mathbb{R};
\end{eqnarray}
\begin{eqnarray}
\phi_2(x_1,\cdots,x_n)=g\left(\prod\limits_{i=1}^n{x_i^{\alpha_i}}\right),\ \alpha_i>0,\ x_i\geq0,
\end{eqnarray}
for  two functions $h$ and $g$. In particular, $\phi_1$ and $\phi_2$ are supermodular functions when $h$ and $g$ are convex functions.  If all $\alpha_i$ are 1, then $\phi_1$ becomes the sum operator and $\phi_2$ becomes the product operator considered in Bignozzi and Puccetti (2015).  By definition, the  distributions $F_1, \cdots, F_n$ which supported by $[0,\infty)$ are $\phi_2$-jointly mixable if and only if the distributions $H_1, \cdots, H_n$ are $\phi_1$-jointly mixable for increasing or decreasing $h$ and $g$, where $H_i(x)=F_i(\exp(\frac{x}{\alpha_i})), 1\le i\le n$; see Lemma 6 in Bignozzi and Puccetti (2015) for the special case of all $\alpha_i$ are 1 and $g(x)=h(x)=x$.

A probability distribution is said to be symmetric if and only if there exists a value $x_0$ such that
$F(x_0-\delta)=F(x_0+\delta)$ for  all $\delta\in {\Bbb R}$.
 That is $X-x_0$ and $-(X-x_0)$ have the same distribution. If $F$ has a density $f$, then $f$ is symmetric if and only if
 $f(x_0-\delta)=f(x_0+\delta)$ for all $\delta\in {\Bbb R}$.
 In particular, if a symmetric distribution has a single mode, the distribution function is a unimodal-symmetric distribution.
 Many commonly used distributions such as normal distribution and $t$-distribution are examples of unimodal-symmetric distributions. According to Cambanis et al.(1981), we know that the class of one-elliptical distributions coincides with the class of one-dimensional symmetric distributions. Therefore, any one-dimensional elliptical distributions which has a single mode is symmetric and unimodal.

Existing results on the mixability of the unimodal-symmetric distributions mainly focused on the complete and joint mixability, which goes back to R\"{u}schendorf and Uckelmann (2002), who gave a simple general construction of random variables in the case of one dimensional unimodal-symmetric distributions, later result on the joint mixability of unimodal-symmetric distribution contains Wang and Wang (2015), who generalized the results of R\"{u}schendorf and Uckelmann (2002) and gave some conditions for the joint mixability of marginal unimodal-symmetric distributions. In this paper, we   consider the $\phi$-joint mixability of unimodal-symmetric distributions.

 In Section 2, we focus on  $\phi_1$-joint  mixability  for the class of   elliptical distributions.  In Section 3, we investigate   $\phi_2$-joint  mixability  for the class of   logarithmic elliptical distributions.

 \section{   $\phi_1$-joint mixability }

In this section, we  consider the $\phi_1$-joint mixability for the class of   elliptical distributions. We first review some  concepts and facts of the  elliptical distributions.

{\bf Definition 2.1} (Fang et al. (1990)). An $n$-dimensional random vector  ${\bf X}$ is said to have an elliptical distribution if its characteristic function can be expressed as
$$E[\exp(i\mathbf{t}'{\bf X})]=\exp(i\mathbf{t}'{\boldsymbol\mu})\psi(\mathbf{t}'\mathbf{\Sigma} \mathbf{t}),\; \mathbf{t}\in \mathbb{R}^n,$$
 where  $\mathbf{{\boldsymbol\mu},\Sigma}$ are parameters, ${\boldsymbol\mu}=(\mu_1,\cdots,\mu_n)'\,\in\,\mathbb{R}^n,$ $\mathbf{\Sigma}=(\sigma_{ij})_{n\,\times\,n}$ is positive semidefinite matrix,  $\psi$:$\,\mathbb{R}\,\to\,\mathbb{R}$ is the characteristic generator of ${\bf X}$. We denote by
 ${\bf X}\sim\text{\Large$\varepsilon$}ll_n(\mathbf{{\boldsymbol\mu},\Sigma},\psi)$.

 Note that not every function $\psi$ can be a characteristic generator, among other things,   it should fulfil the requirement  $\psi(0)=1$.   It is easy to see that if  $\psi(x)=\exp\{-x/2\}$, the elliptical distribution becomes  the normal distribution  $N_n(\mathbf{{\boldsymbol\mu},\Sigma})$.

Let ${\bf \Psi}_n$ be the class of functions $\psi: [0,\infty) \rightarrow \Bbb{R}$
such that function $\psi(|{\bf t}|^2), {\bf t}\in \Bbb{R}^n$  is an $n$-dimensional
characteristic function. It is clear that
${\bf \Psi}_n\subset {\bf \Psi}_{n-1}\cdots\subset {\bf \Psi}_1.$
Denote by ${\bf \Psi}_{\infty}$  the set of characteristic generators that generate an
$n$-dimensional elliptical distribution for arbitrary $n\ge 1$. That is
${\bf \Psi}_{\infty}=\cap_{n=1}^{\infty}{\bf \Psi}_{n}.$ Clearly, if  $\psi(x)=\exp\{-x/2\}$, then $\psi\in {\bf \Psi}_{\infty}$.

\begin{remark} The moments of ${\bf X}\sim\text{\Large$\varepsilon$}ll_n(\mathbf{{\boldsymbol\mu},\Sigma},\psi)$ do not necessarily exist, if $E[X_i]$ exists, it will be given by
$E[X_i]=\mu_i.$
If $Cov[X_i,X_j]$ and/or $Var[X_i]$ exist, they will be given by
$$Cov(X_i,X_j)=-2\psi'(0)\sigma_{ij},$$
and/or
$$Var[X_i]=-2\psi'(0)\sigma^2_i,$$
where $\psi'$ denotes the first derivative of   $\psi$.
 In particular, if the covariance matrix of ${\bf X}$ exists, then it is given by -2$\psi'(0)\mathbf{\Sigma}$. A necessary condition for this covariance matrix to exist is
$|\psi'(0)|<\infty,$
see Cambanis, Huang and Simons (1981).
\end{remark}

For ${\bf X} \sim\text{\Large$\varepsilon$}ll_n(\mathbf{{\boldsymbol\mu},\Sigma},\psi)$, we have the following property.
\begin{lemma} (Fang et al. (1990)).  A random vector  ${\bf X} \sim\text{\Large$\varepsilon$}ll_n(\mathbf{{\boldsymbol\mu},\Sigma},\psi)$ if and only if for any ${\boldsymbol \alpha}=(\alpha_1,\cdots,\alpha_n)'\in\mathbb{R}^n $, ${\boldsymbol \alpha}'{\bf X}\sim\text{\Large$\varepsilon$}ll_1( \mathbf{{\boldsymbol\alpha}}'{\boldsymbol\mu},{\boldsymbol\alpha}'\Sigma{\boldsymbol\alpha},\psi).$
\end{lemma}

For ${\bf X}\sim\text{\Large$\varepsilon$}ll_n(\mathbf{{\boldsymbol\mu},\Sigma},\psi)$,   if the density
exists, then, it has the following form (cf. Fang et al. (1990)):
$$f_{\bf X}({\bf x})=c_n |\mathbf{\Sigma}|^{-\frac12}g_n\left(\mathbf{{(\bf x-\boldsymbol\mu})^T\Sigma}^{-1}{(\bf x-\boldsymbol\mu}\right),$$
for some function $g_n$ called the density generator which satisfies the condition
$$\int_0^{\infty}x^{\frac{n}{2}-1}g_n(x)dx<\infty.$$
 The normalizing constant $c_n$ is given by
\begin{equation*}
c_n=\frac{\Gamma(\frac{n}{2})}{\pi^{\frac{n}{2}}}\left[\int_0^{\infty}x^{\frac{n}{2}-1}g_n(x)dx\right]^{-1}.
\end{equation*}
One may also similarly introduce
the elliptical distribution by the density generator and then write
${\bf X}\sim\text{\Large$\varepsilon$}ll_n(\mathbf{{\boldsymbol\mu},\Sigma},g_n)$.

With the property, we have the following theorem.
\begin{theorem} Assume that $F\sim  \text{\Large$\varepsilon$}ll_1(\mu,\sigma^2,\phi)$, where  $\phi\in{\bf \Psi}_{\infty}$.
Then, $F$  is $\phi_1$-completely mixable with index $n$ for any $n\ge 2$, where $\phi_1$ is defined in (1.2).
\end{theorem}
\emph{Proof}.   Suppose   ${\bf X}\sim\text{\Large$\varepsilon$}ll_n(\mathbf{{\boldsymbol\mu},\Sigma},\phi),   {\boldsymbol\alpha}=(\alpha_1,\cdots,\alpha_n)', \alpha_i\geq0,$
where
${\boldsymbol \mu}=(\mu,\cdots,\mu)'$ and
\begin{eqnarray}
\mathbf{\Sigma}=\sigma^2
       \left(\begin{array}{cccc}
             1&\rho&\ldots&\rho\\
             \rho&1&\ldots&\rho\\
             \vdots&\vdots&\ddots&\vdots\\
             \rho&\rho&\ldots&1\\
       \end{array}\right).
\end{eqnarray}
Here $\rho$, the correlation coefficient of $X_i$ and $X_j$ ($i,j=1,2,\cdots,n)$,  is given as
$$\rho=-\frac{\sum\limits_{i=1}^n{\alpha_i^2}}{2\sum\limits_{1\le i<j\le n}\alpha_i\alpha_j}.$$
   Then $X_i\sim  \text{\Large$\varepsilon$}ll_1(\mu,\sigma^2,\phi)$ ($i=1,2,\cdots n$) and
  $$\sum\limits_{i=1}^n\alpha_i{X_i}\sim    \text{\Large$\varepsilon$}ll_1(\mathbf{{\boldsymbol\alpha}}'{\boldsymbol\mu},{\boldsymbol\alpha}'\Sigma{\boldsymbol\alpha}).$$
 We get
\begin{eqnarray*}
{\boldsymbol\alpha}'\Sigma{\boldsymbol\alpha}&=&\begin{pmatrix} {\alpha_1}, & {\alpha_2}, & {\ldots}, & {\alpha_n} \end{pmatrix}
{\sigma^2}\begin{pmatrix} {1} & {\rho} &  {\ldots} & {\rho}\\{\rho} & {1} &  {\ldots} & {\rho}\\{\vdots} & {\vdots} & {\ddots} & {\vdots}\\{\rho} & {\rho} & {\ldots} & {1}  \end{pmatrix}
\begin{pmatrix} {\alpha_1}\\{\alpha_2}\\{\vdots}\\{\alpha_n} \end{pmatrix}\\
&=&\sigma^2\sum\limits_{i=1}^n\alpha_i^2+2 \sigma^2\sum\limits_{1\le i<j\le n}\alpha_i\alpha_j\rho\\
& =&\sigma^2\sum\limits_{i=1}^n\alpha_i^2+2\sigma^2\sum\limits_{1\le i<j\le n}\alpha_i\alpha_j\left(-\frac{\sum\limits_{i=1}^n{\alpha_i^2}}{2\sum\limits_{1\le i<j\le n}\alpha_i\alpha_j}\right)\\
&=&0,
\end{eqnarray*}
from which we have
\begin{eqnarray*}
E[\exp\{it{\boldsymbol\alpha}'{\bf X}\}]&=&\exp\{i{ t}\mathbf{{\boldsymbol\alpha}'{\boldsymbol\mu}}\}\phi({t}\,\mathbf{{\boldsymbol\alpha}'{\Sigma}{{\boldsymbol\alpha}}}{ t})\\
&=&\exp\{i{ t}\mathbf{{\boldsymbol\alpha}'{\boldsymbol\mu}}\}\phi(0)\\
&=&\exp\{i{ t}\mathbf{{\boldsymbol\alpha}'{\boldsymbol\mu}}\},  t\in \mathbb{R}.
\end{eqnarray*}
This leads to
$$P(\alpha_1X_1+\cdots+\alpha_nX_n={\boldsymbol\alpha}'{\boldsymbol\mu})=1.$$
So that
$$P(\phi_1(X_1,\cdots,X_n)=h({\boldsymbol\alpha}'{\boldsymbol\mu}))=P(h(\alpha_1X_1+\cdots+\alpha_nX_n)=h({\boldsymbol\alpha}'{\boldsymbol\mu}))=1.$$
Thus $F$ is $\phi_1$-completely mixable.

 Suppose that $F_i\sim  \text{\Large$\varepsilon$}ll_1(\mu_i,\sigma_i^2,\psi)$, $i=1,2,\cdots, n$ ($n\ge 3$),  where  $\psi$ is a characteristic generator for an $n$-elliptical
distribution. If $h$ is one-to-one, then it follows from Lemma 2.1 and Theorem 3.7 in Wang and Wang (2016) that
 $(F_1, \cdots, F_n)$  is $\phi_1$-jointly mixable if   and only if
\begin{equation}
\sum_{i=1}^n \alpha_i\sigma_i\ge 2\max\{\alpha_1\sigma_1,\cdots,\alpha_n\sigma_n\},
\end{equation}
where $\phi_1$ is defined by (1.2). If $h$ is not  one-to-one, the condition (2.2) is also a sufficient condition for $\phi_1$-joint  mixability.  Furthermore, if  all means of $F_i$'s   exist, then
the joint center of $(F_1, \cdots, F_n)$ is unique; If  all means of $F_i$'s do not  exist, then
the joint centers of $(F_1, \cdots, F_n)$ are not necessarily unique. For example,  Puccetti, Rigo, Wang and Wang  (2018)     obtain a profound result that for every $n\ge 2$, the set of $n$-centers of the standard Cauchy distribution is
the interval
$$\left[-\frac{\log(n-1)}{\pi}, \frac{\log(n-1)}{\pi}\right].$$
For general  Cauchy distributions with the following probability density functions
$$f(x;\mu,\sigma)=\frac{1}{\pi}\frac{\sigma}{(x-\mu)^2+\sigma^2}, \ -\infty<x<\infty,$$
where $\mu\ge 0,\sigma>0$ are parameters,
the set of $n$-centers  is
the interval
$$\left[-\sigma\frac{\log(n-1)}{\pi}+n\mu, \sigma\frac{\log(n-1)}{\pi}+n\mu\right].$$
It follows that any $C$ in
$$\left[\left(-\sigma\frac{\log(n-1)}{\pi}+n\mu\right)^2, \left(\sigma\frac{\log(n-1)}{\pi}+n\mu\right)^2\right]$$
is the $\phi$-center, where $\phi(x_1,\cdots, x_n)=(x_1+\cdots+x_n)^2$.

Contrastively, for a given supermodular function $\phi$,  the center of a set of $\phi$-jointly mixable distributions might not be unique even when the
distributions have finite mean.  Example 1.1 in Bignozzi and Puccetti (2015) shows that the  center is not unique for $\phi$-jointly mixable discrete distributions having finite means with $\phi(x_1,x_2)=(x_1+x_2)^2$.

Generally, the $\phi_1$-joint mixability of unimodal-symmetric distributions can be characterized by the following theorem, which is a directly consequence of Theorem 3.4 in Wang and Wang (2016).

\begin{theorem} Suppose that $F_1,\cdots,F_n$ have unimodal-symmetric densities and finite means. If $h$ is one-to-one, then $(F_1,\cdots,F_n)$ is $\phi_1$-jointly mixable if and only if
\begin{eqnarray}
\sum\limits_{i=1}^n{\alpha_i}\geq2\max\{{\alpha_1},\cdots,{\alpha_n}\},
\end{eqnarray}
where $\phi_1$ is defined by (2.1).
\end{theorem}

Note that, the center of a set $\phi$-jointly mixable distributions might not be unique even when the distributions have finite first moment. The following theorem gives a sufficient condition for uniqueness of the center of $\phi_1$-joint mixability for continuous unimodal-symmetric distributions.

\begin{theorem}
Suppose $F_1,\cdots,F_n$ are univariate continuous distributions with unimodal-symmetric densities and finite means $\mu_1,\cdots,\mu_n$, respectively. If (2.3) holds, then the $n$-tuple $(F_1,\cdots,F_n)(n\geq3)$ is $\phi_1$-jointly mixable with unique center $K$. Moreover, $K=h\left(\sum\limits_{i=1}^n\alpha_i\mu_i\right)$.
\end{theorem}

The proof is easy.
More generally, according to Proposition 2.3 in Puccetti et al. (2017), we have the following result.
\begin{theorem} Suppose $F_1,\cdots,F_n$ are univariate continuous distributions with unimodal-symmetric densities. Assume that (2.3) holds. If at least $n-2$ of $F_1,\cdots, F_n$ have finite means, without loss of generality, assume that $F_1,\cdots, F_{n-2}$ have finite means, then the $n$-tuple $(F_1,\cdots,F_n)$$(n\ge 3)$ is  $\phi_1$-jointly mixable with unique center $K$. Moreover,
$$K=h\left(\sum_{i=1}^{n-2}\alpha_iE(X_i)+E(\alpha_{n-1} X_{n-1}+\alpha_n X_n) \right),$$
where $(X_1,\cdots,X_n)$ is a joint $\phi_1$-mix for $(F_1, \cdots, F_n)$.
\end{theorem}

In the next theorem, we give a basic proberty for $\phi_1$-joint mixability.
\begin{theorem}
Assume $F_1,\cdots,F_n$ are univariate continuous distributions with unimodal-symmetric densities and finite mean $\mu_1,\cdots,\mu_n$, respectively. Let $\parallel\cdot\parallel$ be any law-determined norm. If the $n$-tuple $(F_1,\cdots,F_n)(n\geq0)$ is $\phi_1$-jointly mixable, $X_i\sim F_i, i=1,\cdots,n$, then the following inequality holds
\begin{eqnarray}
2\max\limits_{1\leq j\leq n}\parallel\alpha_jX_j-\alpha_j\mu_j\parallel\leq\sum\limits_{i=1}^n\parallel\alpha_iX_i-\alpha_i\mu_i\parallel.
\end{eqnarray}
\end{theorem}
The next theorem briefly study a type of convex minimization problem related to the $\phi$-joint mixability.
\begin{theorem}
Suppose $F_1,\cdots,F_n$ are univariate continuous distributions with unimodal-symmetric densities and finite means. Assume that (2.2) hold. Then for any $n\geq3$,
$$
\min\limits_{X_1\sim F_1,\cdots,X_n\sim F_n}E\left(h\left(\sum\limits_{i=1}^n\alpha_iX_i\right)\right)=h\left(\sum\limits_{i=1}^n\alpha_iE[X_i]\right),
$$
where $h$ is strictly convex function.
\end{theorem}
The following theorems concern  joint  mixability and $\phi$-joint  mixability.   Especially we give a sufficient condition for  uniqueness of the center of  $\phi$-joint  mixability for some  elliptical distributions.

\begin{theorem} Suppose that $F_i\sim  \text{\Large$\varepsilon$}ll_1(0,\sigma_i^2,\psi)$ ($i=1,2,\cdots, n$)
  have     densities  of the forms
  \begin{equation}
f_i(x;\sigma_i)=\frac{C_i}{2\sigma_i} f\left(\frac{(x-\nu_i)^2}{\sigma_i^2}\right)+ \frac{C_i}{2\sigma_i}  f\left(\frac{(x+\nu_i)^2}{\sigma_i^2}\right), \ -\infty<x<\infty,
\end{equation}
where $C_i$'s are normalizing constants, $\nu_i\ge 0$, $\sigma_i>0$ are parameters  and  $f$  is a density generator satisfying the condition
$$0<\int_0^{\infty}x^{-\frac12}f(x)dx<\infty.$$
If the elliptical distributions  $G_i  (i=1,2\cdots,n)$ with density generator   $f$     are   jointly mixable,  then  $(F_1, \cdots, F_n)$  is $\phi$-jointly mixable  with center $h(\sum_{i=1}^n \nu_i)$, where
   $\phi(x_1,\cdots,x_n)=h(|x_1+\cdots+x_n|)$ for any function $h$ on $[0,\infty)$.
\end{theorem}
{\bf Proof} \ The  joint mixability of $(G_1, \cdots, G_n)$  implies  that there exist $n$ random variables $Y_1,\cdots,Y_n$ such that $Y_i\sim  \text{\Large$\varepsilon$}ll_1(\nu_i,\sigma_i^2,f)$,  $1\,{\leq}\,i\,{\leq}\,n$, and
 $$Y_1+\cdots+Y_n=\sum_{i=1}^n \nu_i,$$
 and
  there exist $n$ random variables $Z_1,\cdots,Z_n$ such that $Z_i\sim\   \text{\Large$\varepsilon$}ll_1(-\nu_i,\sigma_i^2,f)$,  $1\,{\leq}\,i\,{\leq}\,n$, and
 $$Z_1+\cdots+Z_n=-\sum_{i=1}^n \nu_i.$$
  Define $X_i=\xi Y_i+(1-\xi)Z_i$ ($i=1,2,\cdots, n$), where  $P(\xi=1)=P(\xi=0)=1/2$ and,  random variables $Y_i, Z_i,\xi$ are  independent. Then $X_i$ has the  distribution  $F_i$ and $|X_1+\cdots+X_n|=\sum_{i=1}^n \nu_i$. This shows that $(F_1, \cdots, F_n)$  is $\phi$-jointly mixable with center $h(\sum_{i=1}^n \nu_i)$. This completes the proof of Theorem 2.7. $\hfill\square$

\begin{theorem} Suppose that $F_i\sim  \text{\Large$\varepsilon$}ll_1(0,\sigma_i^2,\psi)$ ($i=1,2,\cdots, n$)
  have     densities  of the forms
  \begin{equation}
f_i(x;\sigma_i)=\frac{C_i}{2\sigma_i} f\left(\frac{(x-\nu_i)^2}{\sigma_i^2}\right)+ \frac{C_i}{2\sigma_i}  f\left(\frac{(x+\nu_i)^2}{\sigma_i^2}\right), \ -\infty<x<\infty,
\end{equation}
where $C_i$'s are normalizing constants, $\nu_i\ge 0, \sigma_i>0$ are parameters  and  $f$  is a density generator satisfying the condition
$$0<\int_0^{\infty}x^{-\frac12}f(x)dx<\infty.$$
Suppose the elliptical distributions  $G_i$'s  with density generator   $f$ are unimodal and (2.1) holds. Then  $(F_1, \cdots, F_n)$  is $\phi$-jointly mixable  with center $h(\sum_{i=1}^n \nu_i)$, where
   $\phi(x_1,\cdots,x_n)=h(|\alpha_1x_1+\cdots+\alpha_nx_n|)$ for any function $h$ on $[0,\infty)$ with constants $\alpha_i$'s.
\end{theorem}
{\bf Proof}\
$X_i \sim   \text{\Large$\varepsilon$}ll_1(0, \sigma^2_i,  \psi)$ implies that $\alpha_i X_i \sim   \text{\Large$\varepsilon$}ll_1(0, \alpha_i^2\sigma^2_i,  \psi)$. Using (2.6) the density of  $\alpha_i X_i$ can be   expressed as
 \begin{equation}
f_i(x; \alpha_i\sigma_i)=\frac{C_i}{2\alpha_i\sigma_i} f\left(\left(\frac{x-\nu_i}{\alpha_i\sigma_i}\right)^2\right)+ \frac{C_i}{2\alpha_i\sigma_i}  f\left(\left(\frac{x+\nu_i}{\alpha_i\sigma_i}\right)^2\right), \ -\infty<x<\infty,
\end{equation}
The unimodality  of $G_i$'s and  condition (2.1) imply that
  there exist $n$ random variables $Y_1,\cdots,Y_n$ such that $Y_i\sim\,G_i,$ $1\,{\leq}\,i\,{\leq}\,n$ and
 $$\alpha_1Y_1+\cdots+\alpha_nY_n=\sum_{i=1}^n \nu_i,$$
 and
  there exist $n$ random variables $Z_1,\cdots,Z_n$ such that $Z_i\sim\,G_i,$ $1\,{\leq}\,i\,{\leq}\,n$ and
 $$\alpha_1Z_1+\cdots+\alpha_nZ_n=-\sum_{i=1}^n \nu_i.$$
  Define $X_i=\xi Y_i+(1-\xi)Z_i$, where  $P(\xi=1)=P(\xi=0)=1/2$ and,  random variables $Y_i, Z_i,\xi$ are  independent. Then $X_i$ has the  density (2.6)
  and $$|\alpha_1X_1+\cdots+\alpha_nX_n|=\sum_{i=1}^n \nu_i,$$ which shows that $(F_1, \cdots, F_n)$  is $\phi$-jointly mixable  with center $h(\sum_{i=1}^n \nu_i)$. This completes the proof of Theorem 2.8. $\hfill\square$

For $X_i\sim  \text{\Large$\varepsilon$}ll_1(0,\sigma_i^2,\psi)$ ($i=1,2,\cdots, n$), if
 $$P(\sum_{i=1}^n X_i=-\mu_n)=P(\sum_{i=1}^n X_i=\mu_n)=\frac12$$
  for some constant $\mu_n>0$, then there exist $2n$ random variables $Y_i$ and $Z_i$ ($i=1,2,\cdots, n$) such that
 $X_i=\xi Y_i+(1-\xi)Z_i$ $(i=1,2,\cdots, n)$ and $P(\sum_{i=1}^n Y_i=-\mu_n)=P(\sum_{i=1}^n Z_i=\mu_n)=1$, where  $P(\xi=1)=P(\xi=0)=1/2$, and   random variables $\{Y_i\}, \{Z_i\},\xi$ are  independent. We conjecture that the density of $X_i$ has the form (2.5).  \\
{\bf Conjecture 1}.\,  Suppose that $F_i\sim  \text{\Large$\varepsilon$}ll_1(0,\sigma_i^2,\psi)$ ($i=1,2,\cdots, n$) have  densities $f_i$, if there exist $n$ random variables $X_1,\cdots, X_n$  such that $X_i\sim F_i$ and $\sum_{i=1}^n X_i$ has the following  two-point distribution
$$G_n(x)=\frac12\delta_{-\mu_n}(x)+\frac12\delta_{\mu_n}(x),\;\; x\in(-\infty,\infty),$$
for some $\mu_n>0$, where $\delta_a(\cdot)$ denotes the point mass at $a$. Then $f_i$  will be of the form
  \begin{equation}
f_i(x;\sigma_i)=\frac{C_i}{2\sigma_i} f\left(\frac{(x-\nu_i)^2}{\sigma_i^2}\right)+ \frac{C_i}{2\sigma_i}  f\left(\frac{(x+\nu_i)^2}{\sigma_i^2}\right), \ -\infty<x<\infty,
\end{equation}
where $C_i$'s are normalizing constants, $\sigma_i>0, \nu_i\ge 0$ are parameters  such that $\sum_{i=1}^n \nu_i=\mu_n$ and  $f$  is  a density generator of  1-dimensional elliptical distribution satisfying the condition
$$0<\int_0^{\infty}x^{-\frac12}f(x)dx<\infty.$$

 \begin{corollary}  Suppose  $F_i\sim  \text{\Large$\varepsilon$}ll_1(0,\sigma_i^2,\psi)$  ($\sigma_i>0$) ($i=1,2,\cdots, n$) have finite means. Then there exist $n$ random variables $X_1,\cdots, X_n$  such that $X_i\sim F_i$ and $\sum_{i=1}^n X_i$ has the following two-point distribution
$$G_n(x)=\frac12\delta_{-\mu_n}(x)+\frac12\delta_{\mu_n}(x),\;\; x\in(-\infty,\infty),$$
for some $\mu_n>0$, if and only if
    $P((\sum_{i=1}^n X_i)^2=\mu_n^2)=1$.
 \end{corollary}
 {\bf Proof}\   The ``only if" part is obvious. Now we prove  the
converse implication.  If there exist $n$ random variables $X_1,\cdots, X_n$  such that $X_i\sim F_i$ and $\sum_{i=1}^n X_i=C$, then $C=0$ since $F_i\sim  \text{\Large$\varepsilon$}ll_1(0,\sigma_i^2,\psi)$ ($i=1,2,\cdots, n$) have finite means. Thus by symmetry of $\sum_{i=1}^n X_i$ and  $P((\sum_{i=1}^n X_i)^2=\mu_n^2)=1$ we have
$$P\left(\sum_{i=1}^n X_i=-\mu_n\right)=P\left(\sum_{i=1}^n X_i=\mu_n\right)=\frac12.$$

\begin{theorem} Suppose that $F_i\sim  \text{\Large$\varepsilon$}ll_1(0,\sigma_i^2,\psi)$  ($\sigma_i>0$) ($i=1,2,\cdots, n$)  with $\psi\in {\bf \Psi}_{\infty}$.  Then, the density $f_i$ of $F_i$ does not have the form of (2.5) with $\nu_i>0$.
\end{theorem}
{\bf Proof}\,  Suppose that $f_i$  has the form (2.5)  with $\nu_i>0$, note that  elliptical distributions   with   $\psi\in {\bf \Psi}_{\infty}$ belonging to the class of scale mixture of the  normal
distributions.  Thus in terms of characteristic functions, (2.5) is equivalent to
\begin{equation}
\int_0^{\infty}\exp\left(-\frac{(\theta\sigma_i)^2}{2}t^2\right)dH_i(\theta)=\cos(t\nu_i)\phi_i\left(\frac{\sigma_i^2 t^2}{2}\right),\ -\infty<t<\infty,
\end{equation}
where $H_i$  a distribution function on $(0, \infty)$  and $\phi_i$ is the characteristic
generator of $g_i$.
 This is a contradiction since the left hand side of (2.6) is  positive for any $t$, but the right hand side of (2.6) is  zero for some $t$.

\begin{remark}    Examples of  elliptical distributions   with   $\psi\in {\bf \Psi}_{\infty}$  are normal distribution,  $T$-distribution, Cauchy distribution, stable laws distribution and double exponential
distribution, and so on. Note that the condition   $\psi\in {\bf \Psi}_{\infty}$  in Theorem 2.3 can be replaced with the condition that  all $F_i$ have positive  characteristic functions.
\end{remark}
We have the following corollary.
\begin{corollary}  Suppose that $F_i\sim  \text{\Large$\varepsilon$}ll_1(0,\sigma_i^2,\psi)$  ($\sigma_i>0$) ($i=1,2,\cdots, n$)   ($\psi\in {\bf \Psi}_{\infty}$)  with finite means.  If there exist $n$ random variables $X_1,\cdots, X_n$  such that $X_i\sim F_i$ and $P((\sum_{i=1}^n X_i)^2=\mu_n^2)=1$ for some constant $\mu_n$,  then  $\mu_n=0$.
\end{corollary}

The following example illustrates that the above situation is not always true.

 {\bf Example 2.1} Suppose  $F$  has  the following  symmetric density
   $$f(x)=\frac{1}{\sqrt{2\pi}}e^{-\frac{x^2+\mu^2}{2}}\cosh(\mu x),\, \ -\infty<x<\infty,$$
   where $\mu>0$ is a constant. Then there exist $n$ random variables $X_1,\cdots, X_n$  such that $X_i\sim F$ and $P((\sum_{i=1}^n X_i)^2=n^2\mu^2)=1$.
 In fact, $f$ can be rewritten as the form
  $$f(x)=\frac{1}{2\sqrt{2\pi}}e^{-\frac{(x+\mu)^2}{2}}+\frac{1}{2\sqrt{2\pi}}e^{-\frac{(x-\mu)^2}{2}},\ -\infty<x<\infty.$$
 Obviously, $f(-x)=f(-x)$, $f$ is unimodal   for  $\mu\le 1$ and bimodal for  $\mu>1$ (cf. Robertson and Fryer (1969)). By Theorem 2.1,
  $F$ is $\phi$-completely mixable with center $(n\mu)^2$, where $\phi(x_1,\cdots,x_n)=(x_1+\cdots+x_n)^2.$
On the other hand, $F$ is the distribution of $\zeta+\eta$, where random variable $\zeta$ and $\eta$ are independent and,  $\zeta\sim N(0,1)$ and  $P(\eta=\mu)=P(\eta=-\mu)=1/2$. Consequently, $0$ is the completely center as well as  the $\phi$-completely mixable center of $F$.

 \section{ $\phi_2$-completely mixability}

 In this section, we will consider the $\phi_2$-completely mixability of  the $\log$-elliptical distributions.
 For any $n$-dimensional vector ${\bf X}=(X_1,\cdots,X_n)'$ with positive components $X_i,$ we define $\log{{\bf X}}=(\log{X_1}, \log{X_2},\cdots,\log{X_n})'.$

\noindent{\bf Definition 3.1}  (Valdez et al. (2009)).  The random vector ${\bf X}$ is said to have a multivariate $\log$-elliptical distribution with
parameters  ${\boldsymbol{\mu}}$ and $\mathbf{\Sigma}$, denoted by ${\bf X} \sim LE_n({\boldsymbol{\mu}},\mathbf{\Sigma},\psi)$, if $\log{\bf X}$ has a multivariate elliptical distribution, i.e.
$\log{\bf X}\sim\text{\Large$\varepsilon$}ll_n({\boldsymbol{\mu}},\mathbf{\Sigma},\psi).$
In particular, when $\psi(x)=\exp\{-x/2\}$, the $\log$-elliptical distributions become  $\log$-normal distributions $LN_n(\mathbf{{\boldsymbol \mu},\Sigma})$.

Using Lemma 2.1 one easy to get
\begin{lemma}  A random vector ${\bf X} \sim LE_n({\boldsymbol{\mu}},\mathbf{\Sigma},\psi)$ if and only if for any  $(\alpha_1,\cdots,\alpha_n)'\in\mathbb{R}^n $, $\prod\limits_{i=1}^n{X_i^{\alpha_i}}\sim LE_1( \mathbf{{\boldsymbol\alpha}}'{\boldsymbol\mu},{\boldsymbol\alpha}'\Sigma{\boldsymbol\alpha},\psi).$ In particular,
any marginal distribution of a $\log$-elliptical
distribution is again $\log$-elliptical.
\end{lemma}

The following theorem shows the $\phi_2$-mixability of $\log$-normal distribution.

\begin{theorem}
Assume that $F$ has the $\log$-normal distribution  $LN_1(\mu,\sigma^2)$. Then $F$ is  $\phi_2$-completely mixable with index $n$ for $n\ge 2$, where $\phi_2$ is defined in (1.3).
\end{theorem}
\emph{Proof}. The theorem is a direct consequence of Theorem 2.1. In fact,  for an $n$-dimensional random vector ${\bf X}\sim LN_n(\mathbf{{\boldsymbol \mu},\Sigma})$,
where ${\boldsymbol\mu}$ and ${\boldsymbol\Sigma}$ are the same as in Theorem 2.1, then $X_i\sim LN_1(\mu,\sigma^2 )$ $(i=1,2,\cdots,n)$  and $\ln{{\bf X}}{\sim}N_n(\mathbf{{\boldsymbol \mu},\Sigma})$.   According to Theorem 2.1, we have
$$P\left(\sum\limits_{i=1}^n\alpha_i\ln{X_i}={\boldsymbol\alpha}'{\boldsymbol\mu}\right)=1,$$
and hence
$$P\left(\prod\limits_{i=1}^n{X_i^{\alpha_i}}=e^{{\boldsymbol{\alpha}}'{\boldsymbol{\mu}}}\right)=1.$$
It follows that
$$P\left[g\left(\prod\limits_{i=1}^n{X_i^{\alpha_i}}\right)=g(e^{{\boldsymbol{\alpha}}'{\boldsymbol{\mu}}})\right]=1.$$
That is
$$P[\phi_2(X_1,\cdots,X_n)=g(e^{{\boldsymbol{\alpha}}'{\boldsymbol{\mu}}})]=1,$$
and thus $F$ is $\phi_2$-completely mixable.

 In the following we consider the $\log$-elliptical distributions, which can be seen as direct extensions  of  log-normal distributions.

\begin{theorem}
The $\log$-elliptical distribution $F\sim LE_1(\mu,\sigma^2,\phi)$   is $\phi_2$-completely mixable with index $n$ for $n\ge 2$, where    $\phi$  is a characteristic generator for an $n$-elliptical distribution ($n\ge 2$) and  $\phi_2$ is defined in (1.3).
\end{theorem}
\emph{Proof}. Using the same argument as that in the proof to Theorem 3.1, the result directly follows from   Theorem 2.2.

Using the result of the last section, one finds that if $g$ is one-to-one, then
 $(F_1, \cdots, F_n)$  is $\phi_2$-jointly mixable if   and only if
$$\sum_{i=1}^n \alpha_i\sigma_i\ge 2\max\{\alpha_1\sigma_1,\cdots,\alpha_n\sigma_n\},$$
where $\phi_2$ is defined by (1.3).
 Furthermore, if  all means of $F_i$'s   exist, then
the $\phi_2$-jointly   center of $(F_1, \cdots, F_n)$ is unique; If  all means of $F_i$'s do not  exist, then
the  $\phi_2$-jointly  centers of $(F_1, \cdots, F_n)$ are not necessarily unique. For example,
using the result of Puccetti, Rigo, Wang and Wang  (2018) in the last section we have that   for every $n\ge 2$, the set of  $n$-product centers of  the log-Cauchy distribution with the probability density function
 $$f(x;\mu,\sigma)=\frac{1}{x\pi}\frac{\sigma}{(\log x-\mu)^2+\sigma^2}, \ x>0,$$
 is the interval
$$\left[\exp\left(-\sigma\frac{\log(n-1)}{\pi}+n\mu\right), \exp\left(\sigma\frac{\log(n-1)}{\pi}+n\mu\right)\right],$$
where $\mu\ge 0$   and $\sigma>0$ are parameters.

Corresponding to Theorems 2.1 and 2.2, we have the following theorems.
\begin{theorem} Suppose that $F_i\sim  LE_1(0,\sigma_i^2,\psi)$ ($i=1,2,\cdots, n$)
  have     densities  of the forms
  \begin{equation}
f_i(x;\sigma_i)=\frac{C_i}{2\sigma_i x} f\left(\frac{(\log x-\nu_i)^2}{\sigma_i^2}\right)+ \frac{C_i}{2\sigma_i x}  f\left(\frac{(\log x+\nu_i)^2}{\sigma_i^2}\right), \ x>0,
\end{equation}
where $C_i$'s are normalizing constants, $\nu_i\ge 0$, $\sigma_i>0$ are parameters  and  $f$  is a density generator satisfying the condition
$$0<\int_0^{\infty}x^{-\frac12}f(x)dx<\infty.$$
If the $\log$-elliptical distributions  $G_1,\cdots, G_n$ with  density generator  $f$  are jointly mixable,  then  $(F_1, \cdots, F_n)$  is $\phi$-jointly mixable  with center $g(\sum_{i=1}^n \nu_i)$, where
   $\phi(x_1,\cdots,x_n)=g(|\log(\prod_{i=1}^n x_i)|)$ for any function $g$ on $[0,\infty)$.
\end{theorem}

\begin{theorem}  Suppose that $F_i\sim  LE_1(0,\sigma_i^2,\psi)$ ($i=1,2,\cdots, n$)
  have     densities  of the forms
  \begin{equation}
f_i(x;\sigma_i)=\frac{C_i}{2\sigma_i x} f\left(\frac{(\log x-\nu_i)^2}{\sigma_i^2}\right)+ \frac{C_i}{2\sigma_i x}  f\left(\frac{(\log x+\nu_i)^2}{\sigma_i^2}\right), \ x>0,
\end{equation}
where $C_i$'s are normalizing constants, $\nu_i\ge 0$, $\sigma_i>0$ are parameters  and  $f$  is  a density generator satisfying the condition
$$0<\int_0^{\infty}x^{-\frac12}f(x)dx<\infty.$$
Suppose the $\log$-elliptical distributions  $G_i$'s  with density generator   $f$ are unimodal and (2.1) holds. Then  $(F_1, \cdots, F_n)$  is $\phi$-jointly mixable  with center $g(\sum_{i=1}^n \nu_i)$, where
   $$\phi(x_1,\cdots,x_n)=g\left(\left|\log\left(\prod_{i=1}^n x_i^{\alpha_i}\right)\right|\right)$$
    for any function $g$ on $[0,\infty)$ and constants $\alpha_i>0$.
\end{theorem}
We next give a lemma, which will be crucial to construct conditions on $\phi_2$-joint mixability.

\begin{lemma}  Suppose that $F_1,\cdots,F_n$ are $n$ univariate continuous distribution functions, $g$ is a nonconstant continuous function on $\mathbb{R}^n$, $f$ is a nonconstant continuous function on the real line. Then there exist n random variables $X_i\sim F_i$, $i=1,\cdots,n$ and constant $C$ such that
\begin{eqnarray*}
P(f(g(X_1,\cdots,X_n)))=f(C))=1,
\end{eqnarray*}
if and only if
\begin{eqnarray*}
P(g(X_1,\cdots,X_n))=C)=1,
\end{eqnarray*}
\end{lemma}

By using Lemma 3.1, we obtain the following theorem.

\begin{theorem} Suppose $F_1,\cdots,F_n$ are univariate distributions with unimodal-symmetric densities and finite means, $H_i(x)=F_i(\log x)(i=1,\cdots,n)$ are positive distributions. If $g$ is one-to-one, then $(H_1,\cdots,H_n)$ is $\phi_2$-jointly mixable if and only if
\begin{eqnarray}
\sum\limits_{i=1}^n{\alpha_i}\geq2\max\{{\alpha_1},\cdots,{\alpha_n}\},
\end{eqnarray}
where $\phi_2$ is defined by (3.1).
\end{theorem}

\emph{Proof}. Suppose random vector ${\bf X}=(X_1,\cdots,X_n)',X_i\sim H_i$, $\boldsymbol{\alpha}=(\alpha_1,\cdots,\alpha_n)'$, $\alpha_i>0$. According to Lemma 3.1, there exist constant $C>0$ and $n$ random variables $X_1,\cdots,X_n$ such that
$$
P(\prod\limits_{i=1}^nX_i^{\alpha_i}=C)=1,
$$
if and only if
$$
P(f(\prod\limits_{i=1}^nX_i^{\alpha_i})=f(C))=1,
$$
for any noncontinuous function $f$ on real line.
Setting $f(x)=\log x$, we have that
$$
P(\sum\limits_{i=1}^n{\alpha_i}\log X_i=\log C)=1,
$$
by Theorem 3.1, it holds if and only if
$$
\sum\limits_{i=1}^n{\alpha_i}\geq2\max\{{\alpha_1},\cdots,{\alpha_n}\}.
$$
The result follows since $g$ is one-to-one.

As mentioned earlier, the center of a set $\phi$-jointly mixable distributions might not be unique, similar to Theorem 2.2, we have the following result.

\begin{theorem} Suppose $F_1,\cdots,F_n$ are continuous distributions with unimodal-symmetric densities and finite means $\mu_1,\cdots,\mu_n$, respectively, $H_i(x)=F_i(\log x)$ $(i=1,\cdots,n)$ are positive distributions. Assume that (3.2) holds. Then the $n$-tuple $(H_1,\cdots,H_n)(n\geq3)$ is $\phi_2$-jointly mixable with unique center $K$. Moreover, $K=g(\prod\limits_{i=1}^n \exp(\alpha_i\mu_i))$.
\end{theorem}

\emph{Proof}. Suppose random vector ${\bf X}=(X_1,\cdots,X_n)', X_i\sim H_i$. From the proof of Theorem 3.1, we see that there exist constant $C$ and $n$ random variables $X_i(i=1,\cdots,n)$ such that
$$
P(\prod\limits_{i=1}^nX_i^{\alpha_i}=C)=1.
$$
Moreover,
$$
C=\exp\left(\sum\limits_{i=1}^n\alpha_iE(\log X_i)\right)=\exp\left(\sum\limits_{i=1}^n\alpha_i\mu_i\right),
$$
if and only if
$$
\sum\limits_{i=1}^n{\alpha_i}\geq2\max\{{\alpha_1},\cdots,{\alpha_n}\},
$$
Therefore,
\begin{eqnarray}
P\left(g\left(\prod\limits_{i=1}^nX_i^{\alpha_i}\right)=K\right)=1,
\end{eqnarray}
where $K=g(C)$.
Since $F_1,\cdots,F_n$ are continuous, so that $H_i(x)(i=1,\cdots,n)$ are continuous for any $n$ random variables $X_i\sim H_i$, $i=1,\cdots,n$. Then $\prod\limits_{i=1}^nX_i^{\alpha_i}$ is either degenerate or is a continuous random variable. Assume $Var\left(\prod\limits_{i=1}^nX_i^{\alpha_i}\right)<\infty$, whenever
$$
Var\left(\prod\limits_{i=1}^nX_i^{\alpha_i}\right)=0,
$$
we have found $n$ random variable $X_1,\cdots,X_n$, $X_i\sim H_i(i=1,\cdots,n)$ such that (3.3) holds.
Whenever
$$
Var\left(\prod\limits_{i=1}^nX_i^{\alpha_i}\right)>0,
$$
$\prod\limits_{i=1}^nX_i^{\alpha_i}$ is a continuous random variable which supported on $I$. Thus for any piecewise continuous function $g$ on $I$, $g\left(\prod\limits_{i=1}^nX_i^{\alpha_i}\right)$ is also a continuous random variable. It follows that for any constant $K$,
$$
P\left(g\left(\prod\limits_{i=1}^nX_i^{\alpha_i}\right)=K\right)=0.
$$
This ends the proof.

Corresponding to Theorems 2.3-2.5, we have the following results.

\begin{theorem} Suppose $F_1,\cdots,F_n$ are continuous distributions with unimodal-symmetric densities, $H_i(x)=F_i(\log x)(i=1,\cdots,n)$ are positive distributions. Assume that (3.2) holds. If at least $n-2$ of $F_1,\cdots, F_n$ have finite means, without loss of generality, assume that $F_1,\cdots,F_{n-2}$ have finite means, then the $n$-tuple $(H_1,\cdots,H_n)$$(n\ge 3)$ is $\phi_2$-jointly mixable with unique center $K$. Moreover,
$$K=g\left(\exp\left(\sum_{i=1}^{n-2}\alpha_iE(\log X_i)+E(\alpha_{n-1} \log X_{n-1}+\alpha_n \log X_n)\right) \right),$$
where $(X_1,\cdots,X_n)$ is a joint $\phi_2$-mix for $(H_1, \cdots, H_n)$.
\end{theorem}

\begin{theorem} Assume the univariate distributions $F_i(i=1,\cdots,n)$ are continuous distributions with unimodal-symmetric densities and finite means $\mu_1,\cdots,\mu_n$, respectively, $H_i(x)=F_i(\log x)(i=1,\cdots,n)$ are positive distributions. Let $\parallel\cdot\parallel$ be any law-determined norm. If the $n$-tuple $(H_1,\cdots,H_n)(n\geq0)$ is $\phi_2$-jointly mixable, $X_i\sim H_i$, $i=1,\cdots,n$, then the following inequality holds
\begin{eqnarray}
2\max\limits_{1\leq j\leq n}\parallel\alpha_j\log X_j-\alpha_j\mu_j\parallel\leq\sum\limits_{i=1}^n\parallel\alpha_i\log X_i-\alpha_i\mu_i\parallel.
\end{eqnarray}
\end{theorem}

{\bf Proof}. From the proof of Theorem 3.2, we see that the $n$ tuple $(H_1,\cdots,H_n)$ is $\phi_2$-jointly mixable, if and only if $(H_1,\cdots,H_n)$ is $\phi_2^1$-jointly mixable, where
$\phi_2^1(x_1,\cdots,x_n)=\prod\limits_{i=1}^nx_i^{\alpha_i}$, that is there exist random variables $Y_1\sim H_1,\cdots,Y_n\sim H_n$ such that
$$
P\left(\prod\limits_{i=1}^nY_i^{\alpha_i}=C\right)=1,
$$
where $C=\prod\limits_{i=1}^n\exp(\alpha_i\mu_i)$.
Thus
$$
P\left(\sum\limits_{i=1}^n{\alpha_i}\log Y_i=\sum\limits_{i=1}^n\alpha_i\mu_i\right)=1,
$$
that is
$$
P\left(\sum\limits_{i=1}^n\Bigg(\alpha_i\log Y_i-\alpha_i\mu_i\Bigg)=0\right)=1,
$$
from which we have
$$||\alpha_j \log Y_j-\alpha_j\mu_j||\le \sum_{i=1, i\neq j}^n ||\alpha_i\log Y_i-\alpha_i\mu_i||.$$
Thus (3.4) follows since $||\cdot||$ is law determined.

\begin{theorem} Suppose that the univariate distributions $F_i(i=1,\cdots,n)$ are continuous distributions with unimodal-symmetric densities and finite means $\mu_1,\cdots,\mu_n$, respectively, $H_i(x)=F_i(\log x)(i=1,\cdots,n)$ are positive distributions. Assume that (3.2) holds. Then for any $n\geq3$,
$$
\min\limits_{X_1\sim F_1,\cdots,X_n\sim F_n}E\left(g\left(\prod\limits_{i=1}^nX_i^{\alpha_i}\right)\right)=g\left(\exp\left(\sum\limits_{i=1}^n\alpha_i\mu_i\right)\right),
$$
where $g$ is piecewise continuous function.
\end{theorem}

\noindent{\bf Data Availability}: There are no data used in this study.

\noindent{\bf Author Contributions}: These authors contributed equally to this work.

\noindent{\bf Conflicts of Interest}: The authors declare that there are no conflicts of interest regarding the publication of this article.

\noindent{\bf Acknowledgements}. %The authors wish to acknowledge the comments and suggestions made by the anonymous referees   which helped in improving this version of the paper.
The research was supported by the National Natural Science Foundation of China (No.  12071251).

%\clearpage

\end{document}